\documentclass[a4paper, leqno, 12pt]{article}
\usepackage{amsmath, amssymb}
\usepackage{graphicx}
\numberwithin{equation}{section}
\begin{document}
\begin{center}
{\Large
Phase Field Equation in the Singular Limit of Stefan Problem
}
\end{center}
\begin{center} 
Jun-ichi Koga, Jiro Koga and Shunji Homma
\end{center}
\begin{center}
Division of Materials Science, Graduate School of Science and Engineering,
Saitama University, Sakura-ku, Saitama, 338-8570 Japan
\end{center}
\begin{abstract}
\qquad The classical Stefan problem is reduced as the singular limit of phase-field equations. These equations are for temperature $u$ and the phase-field $\varphi$, consists of a heat equation:
$$
u_t+\ell\varphi_t=\Delta u,
$$
and a Ginzburg-Landau equation:
$$
\epsilon\varphi_t=\epsilon\Delta\varphi -\frac{1}{\epsilon}W^\prime (\varphi )+\ell (\varphi )u,
$$
where $\ell$ is a latent heat and $W$ is a double-well potential whose wells, of equal depth, correspond to the solid and liquid phases. \\
\qquad When $\epsilon\to 0$, the velocity of the moving boundary $v$ in one dimension and  that of the radius in the cylinder or sphere are shown as the following Stefan problem,\\
$$
\left\{
\begin{array}{l}
u_t-\Delta u =0\\\\
\displaystyle v=\frac{1}{2}\left[\frac{\partial u}{\partial n}\right]_\Gamma\\\\
\displaystyle u=-\frac{m}{2\ell}[\kappa -\alpha v]_\Gamma
\end{array}
\right.
$$
where $\alpha$ is a positive parameter, $[\frac{\partial u}{\partial n}]_\Gamma$ is the jump of the normal derivatives of $u$ (from solid to liquid), and $m=\int_{-1}^1\left(2W(\varphi)\right)^{1/2}d\varphi$.\\
\qquad Since it is sufficient to describe the phase transition of single component by the phase-field equation, we analyze the phase-field equation,
$$
\frac{\partial\varphi}{\partial t}=a^2\Delta\varphi+f\left(\varphi\right),
$$
where $a$ is a positive parameter, and $f(\varphi)$ a function of the double-well potential and investigate whether the equation shows the Stefan problem or not. The velocity of the moving boundary in the cylinder and sphere are determined and the result of the simulation of the equation is also presented.  Next, we consider the velocity of interface which depends on the temperature. The control of width of diffusion layer by the parameter of phase-field equations is investigated in order to realize the singular limit of phase field equations by the numerical method. 
\end{abstract}
\section{Introduction}
Phase-field equations are used in the field of phase transition.  Recently, the relationships between the physical properties and the parameters of these equations have been determined and as a result, phase-field equations have been established as the position of the model for simulation of the phase transition behavior in engineering field.\\
\qquad Phase-field equations are also applied to the double layer model of interface between fluids for example, interface of liquid-gas and liquid-liquid contact.  However in the cases of contact of gas-sold and liquid-solid, we should modify these equations for the diffusion layer not to be double, but the width of diffusion layer tends to zero.\\
\qquad The singular limit of the phase-field equations have been investigated theoretically and numerically.  Caginalp and a co-researcher [1-4] have studied the reduction of these equations in the singular limit of Stefan problems and Hele-Shaw problem. Furthermore H. Soner [11] showed that the solution of the phase-field equations converged to that of the Stefan problem. An asymptotic behavior of solutions for the Stefan problem is calculated by Carslaw and Jaeger [5], Friedman [6]. In their study, the moving interface distance on the process of the solidification is a half-square of time, namely, this result is coincided with the dependence on time as the free moving interface length, with two points from the origin, of the phase-field equations. \\
\qquad It is sufficient to describe the phase transition of single component by the phase-field equation. Therefore we analyze the equation and obtain the asymptotic behavior of the equation in order to investigate that the  equation shows the Stefan problem.\\
\qquad This paper consists as follows. Sec.2 is devoted to the review of studies by many researchers in order to  make clear the relations among the theories concerning of this problem for applied mathematician, physicist and engineer. Sec.3 is a numerical example.\\
\section{Review on Mathematics of This Problem}
Caginalp and co-researchers reported that the phase-field equations converged to the equations for the Stefan problem as singular limit.  These equations are for temperature $u$ and the phase-field $\varphi$, consists a heat equation:
$$
u_t+\ell\varphi_t=\Delta u,
$$
and a Ginzburg-Landau equation:
$$
\epsilon\varphi_t=\epsilon\Delta\varphi -\frac{1}{\epsilon}W^\prime (\varphi )+\ell (\varphi )u,
$$
where $\ell$ is a latent heat and $W$ is a double-well potential whose wells, of equal depth, correspond to the solid and liquid phases. \\
\qquad When $\epsilon\to 0$, the velocity of the moving boundary $v$ in one dimension and  that of the radius in the cylinder or sphere are shown as the following Stefan problem and let $\Gamma (t)$ be the interface separating the two regions $\Omega(t) = \{x: \varphi (t, x) = - 1\}$ or $\{\varphi =1\}$,\\
$$
\left\{
\begin{array}{l}
u_t-\Delta u =0,\,\, \mathrm{in}\,\,\mathrm{\Omega = \{x: \varphi (t, x) = \pm 1\}}\\\\
\displaystyle v=\frac{1}{2}\left[\frac{\partial u}{\partial n}\right]_\Gamma\\\\
\displaystyle u=-\frac{m}{2\ell}[\kappa -\alpha v]_\Gamma
\end{array}
\right.
$$ 
where $\alpha$ is a positive parameter, $[\frac{\partial u}{\partial n}]_\Gamma$ is the jump of the normal derivatives of $u$ (from solid to liquid), and $m=\int_{-1}^1\left(2W(\varphi)\right)^{1/2}d\varphi$.  On the other hand, Iida [7] suggested that those equations did to the equations for the two phase Stefan problem by the same singular limit. These of the fundamental idea is based on the work about the model of ``Caginalp type''.  More precisely, these results can be derived by the method of matched asymptotic expansion between the inner expansion and the outer expansion. \\
\qquad Soner reported that the solution of the phase-field equations converges to that of the Stefan problem, when $\epsilon\to 0$. Now, when $\epsilon\to 0$, the solution $(u^\epsilon, \varphi^\epsilon)$ of
$$
\left\{
\begin{array}{l}
u^\epsilon_t+\ell\varphi^\epsilon_t=\Delta u^\epsilon,\\\\
\displaystyle\epsilon\varphi^\epsilon_t=\epsilon\Delta\varphi^\epsilon -\frac{1}{\epsilon}W^\prime (\varphi^\epsilon )+\ell (\varphi^\epsilon )u^\epsilon,
\end{array}
\right.
$$
converges to the solution $(u,\varphi)$ of
$$
\left\{
\begin{array}{l}
u_t-\Delta u= -(h(\varphi))_t\\\\
\displaystyle u=\pm\frac{m}{2\ell}[\kappa -\alpha v]_\Gamma
\end{array}
\right.
$$
where $+$ stands for solidification and $-$ melting, and $\ell (\varphi)=h^\prime (\varphi)=\sqrt{2W(\varphi)}$.\\
\qquad From these results, the moving velocity of interface is obtained from the limit of the phase-field equations. The velocity of interface without the interface tension is 
$$
v=\pm\frac{2\ell}{\alpha m}u,
$$
when the temperature $u=0$ on $\Gamma$, the velocity
$$
v=\pm\frac{1}{\alpha}\kappa
$$
where $\alpha$ is a positive constant.
\section{Numerical Example}
We consider the following equations,
\begin{equation}
\Phi_t(r,t)=\epsilon^2\Delta\varphi +f(\varphi).
\end{equation}
Suppose that the width of the interface between solid and liquid $\epsilon$ is sufficiently small in comparison with the radius of solid $R(t)$ as shown in Fig.1. Let us consider the time variation $R(t)$. In the case that $\epsilon/R\ll 1$, we can write
\begin{equation}
\Phi (r,t)=\Phi_0(r\pm R(t)).
\end{equation}
Substituting (3.2) to (3.1), we get
\begin{equation}
\pm\dot{R}\Phi_0^\prime =\epsilon^2\left(\Phi_0^{\prime\prime}+\frac{d-1}{r}\Phi_0^\prime\right)+f(\Phi_0).
\end{equation}
where $d$ is a dimensional number, 2 or 3.\\
\qquad Since the equilibrium is assumed at the interface, the following equation is given, namely,
\begin{equation}
\epsilon^2 \frac{d^2}{dx^2}\Phi_0+f(\Phi_0)=0
\end{equation}
holds. Here $\Phi_0^\prime$ is a large value at the neighborhood of $r=R$ and rapidly decrease to zero out of that neighbor.   Hence we replace $\Phi_0^\prime/r$ to $\Phi_0^\prime/R$ and as a result, from $(3.3)$, we get
\begin{equation}
\dot{R}\cong\epsilon^2\frac{d-1}{R},
\end{equation}
and
\begin{equation}
(R(t))^2=R_0^2+2\epsilon^2(d-1)t
\end{equation}
where $R_0$ is the initial radius. The relation of $(R_0^2$ vs. $t$ is good agreement with tendency of the exact solutions [5] of the Stefan problem.\\
\begin{figure}[h]
\begin{center}
\includegraphics[width=7cm,clip]
{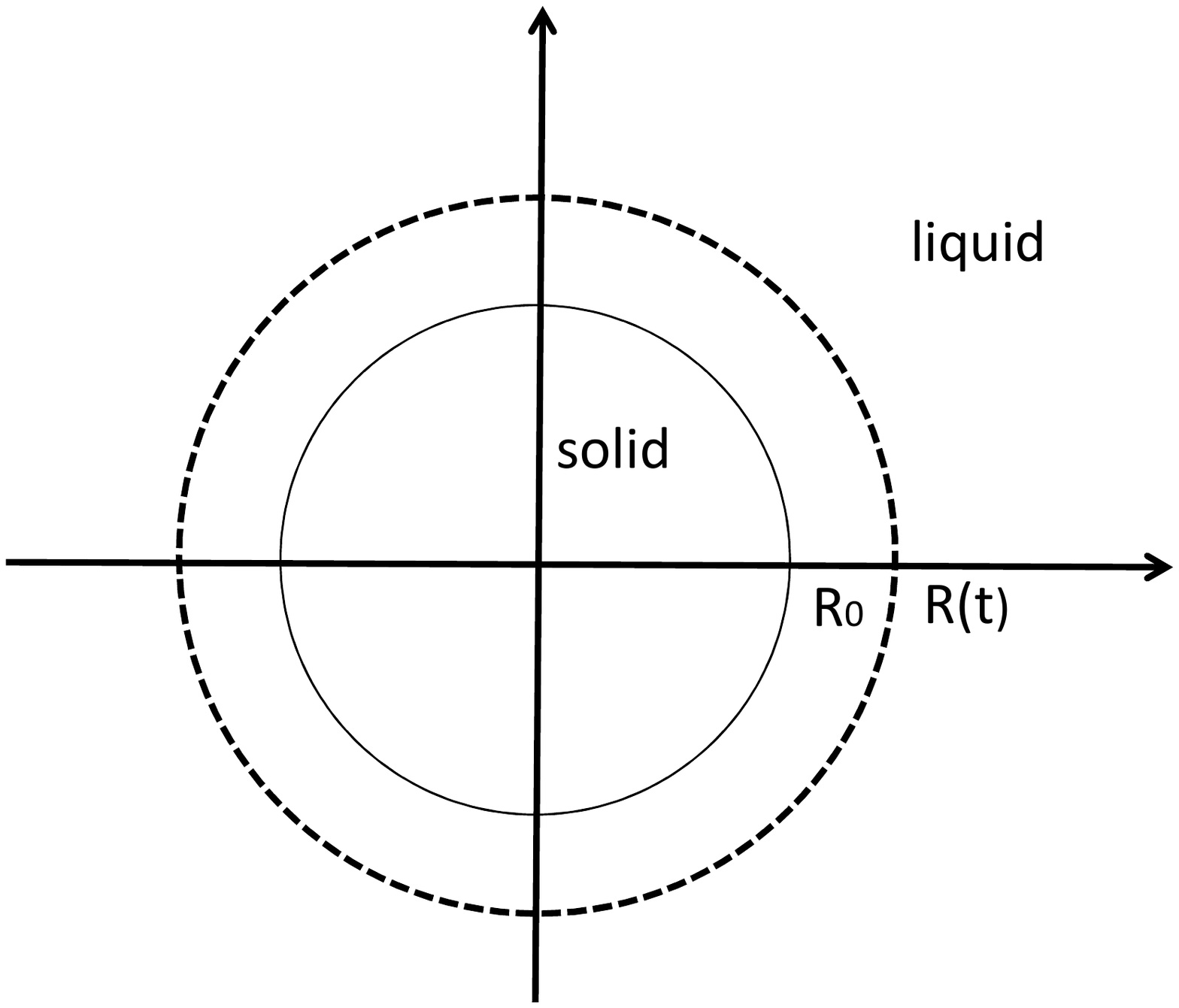}
\caption{$R(t)$, $R_0$ and $x,y$-coordinate}
\end{center}
\end{figure}
\qquad Next, we consider the velocity of interface which depends on the temperature. The control of width of diffusion layer by the parameter of phase-field equations is investigated in order to realize the singular limit of phase field equations by the numerical method. Karma and Rappel [8] in 1996 derived parameters for the so-called thin-interface limit of phase-field equations, where the interface thickness can be controlled to be small, and the classical interface conditions are satisfied for a finite thickness. Their analysis allowed for the first time fully resolved computations to be made for three dimensional dendrites with arbitrary interface kinetics, by Karma and Rappel [9] in 1998.\\
\qquad We used their model to simulate the Stefan problem by the phase-field equations. In this simulation, we gave the function as
\begin{equation}
W(\varphi)=-\frac{\varphi^2}{2}+\frac{\varphi^4}{4}
\end{equation} 
and
\begin{equation}
\frac{\partial u}{\partial t}=\Delta u+\frac{1}{2}\frac{\partial\varphi}{\partial t}
\end{equation}
\begin{equation}
\frac{\partial\varphi}{\partial t}=\Delta\varphi +W^\prime (\varphi)-\lambda (u)
\end{equation}
The right hand side of the equations is introduced to the parameter $\lambda$. The parameter $\lambda$ is controlled by the thickness of diffusion layer between solid and liquid, and is similar to the one used by Caginalp and others (e.g. see MacFadden, {\it et. al.} (1993) [10]).\\
\qquad Now, we consider the steady state of one dimensional phase-field equation to compare the sharp interface model. The boundary conditions of the equations are 
\begin{equation}
V=\frac{\partial u}{\partial x},
\end{equation}
\begin{equation}
u_i=-\beta V,
\end{equation}
\begin{equation}
u(+\infty)=-\delta.
\end{equation}
where in the equation $(3.11)$, subscript $i$ stands for the interface. Using the equations, we obtained the solution of the steady state equations as
\begin{equation}
V=\frac{\delta -1}{\beta}
\end{equation}
\begin{equation}
u=\exp [-Vx]-\delta
\end{equation} 
in the liquid $(x\geqq 0)$ and $u=1-\delta$ in the solid $(x\leqq 0)$.\\
\qquad In one dimension, phase-field equations take the form
\begin{equation}
\frac{\partial\varphi}{\partial t}=\frac{\partial^2 \varphi}{\partial x^2}+\left(\varphi -\varphi^3\right)+\lambda u
\end{equation}
\begin{equation}
\frac{\partial u}{\partial t}=\frac{\partial^2 u}{\partial x^2}+\frac{1}{2}\frac{\partial\varphi}{\partial t}.
\end{equation}
The steady-state growth equations in the moving frame of interface, which yield
\begin{equation}
V\frac{\partial\varphi}{\partial x}+\frac{\partial^2\varphi}{\partial x^2}+\varphi-\varphi^3+\lambda u=0,
\end{equation}
\begin{equation}
V\frac{\partial u}{\partial x}+\frac{\partial^2 u}{\partial x^2}-\frac{V}{2}\frac{\partial\varphi}{\partial x}=0.
\end{equation}
The solution of these equations with $u$ subject to the far field boundary conditions determines the planer interface velocity as a function of undercooling.\\
\qquad We numerically analyzed the equations $(3.17)$, $(3.18)$. Figure 2 and Figure 3 show the profiles of $\varphi$ and $U$. The thickness of the interface is dependent on the parameter $\lambda$. The thickness of $\varphi$ and $U$ with $\lambda$ decreasing. This procedure of thin interface limit gives the solution of the Stefan problem as singular limit of phase-field equations. 
\begin{figure}[h]
\begin{center}
\includegraphics[width=10cm, height=13cm, clip]
{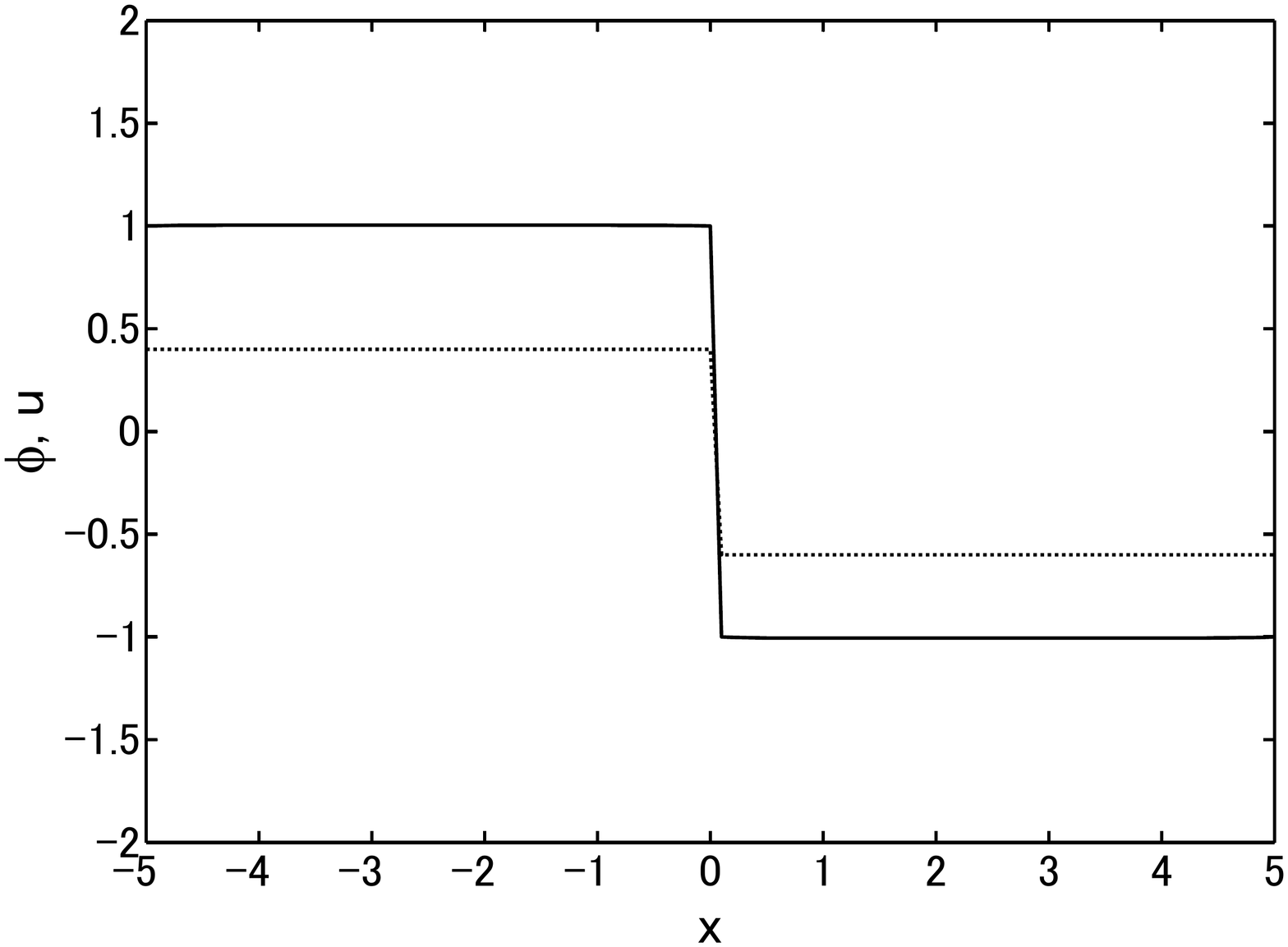}
\caption{Profile of $\phi$, $u$: $\phi:-$, $u:\cdots$, on the conditions that $\beta =0.2572$, $\delta=0.6$, $\lambda =0.3$.}
\end{center}
\end{figure}
\begin{figure}[h]
\begin{center}
\includegraphics[width=10cm, height=13cm, clip]
{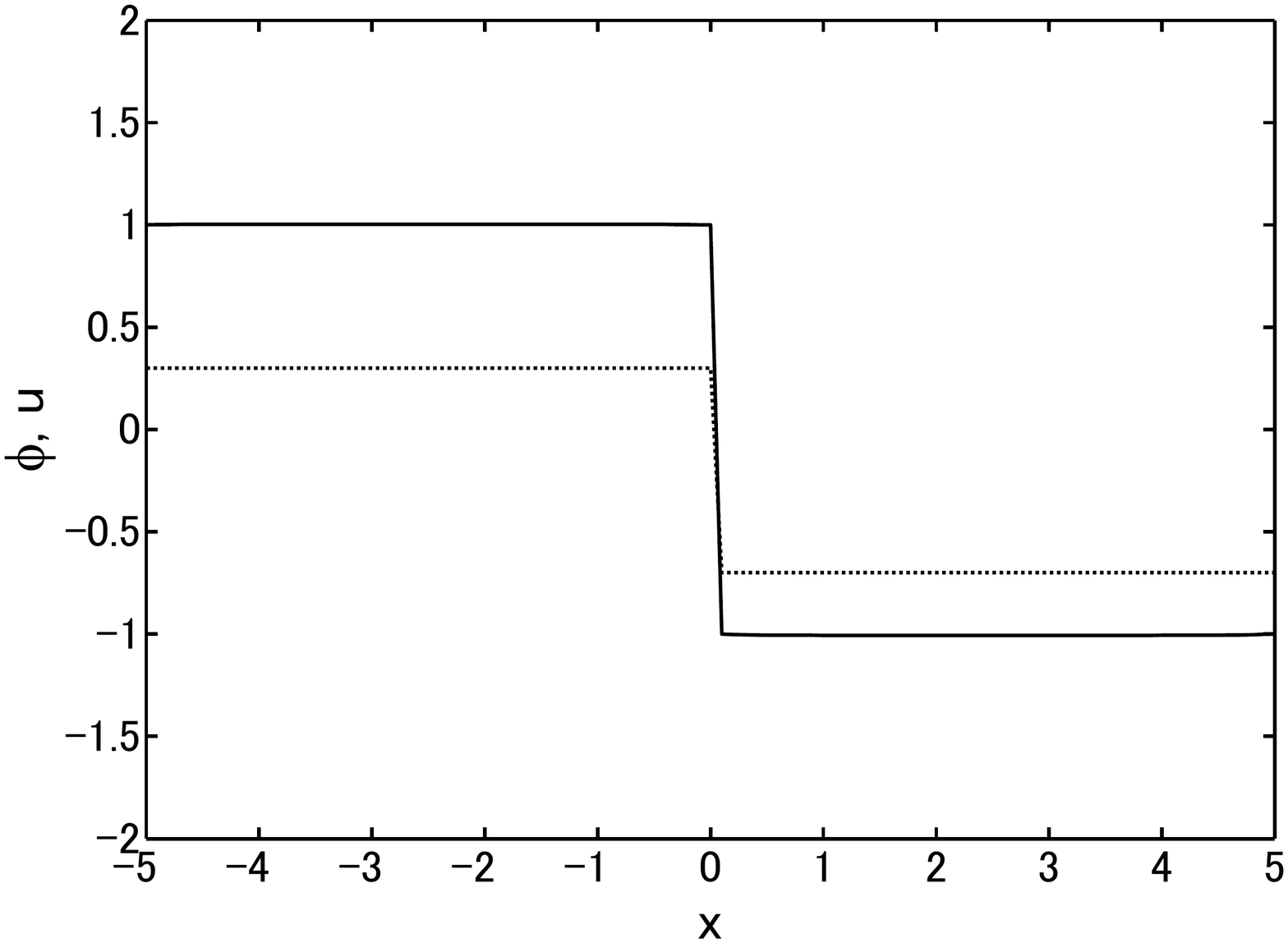}
\caption{Profile of $\phi$, $u$: $\phi:-$, $u:\cdots$, on the conditions that $\beta =0.2572$, $\delta =0.7$, $\lambda =0.1$}
\end{center}
\end{figure}
\section{Conclusion}
\qquad Rearranged the theorems of phase-field equations in singular limit of Stefan problem, we reviewed the previous work reported by mathematician. The velocity of the moving interface between solid and liquid was simply determined.

\end{document}